\documentclass{article}[12pt]
\usepackage[utf8]{inputenc}

\usepackage[utf8]{inputenc}
\usepackage{amsfonts,amsmath,latexsym,amssymb,verbatim,amsbsy,amsthm}
\usepackage{comment}

\usepackage{graphicx}
\usepackage{nicefrac}
\usepackage{dsfont}
\usepackage{mathrsfs}
\usepackage{cancel}

\usepackage[normalem]{ulem}
\usepackage{color}
\usepackage{caption,subcaption}

\usepackage[dvipsnames]{xcolor}
\usepackage{relsize}

\numberwithin{equation}{section}
\numberwithin{figure}{section}
\numberwithin{table}{section}

\theoremstyle{plain}

\theoremstyle{definition}

\theoremstyle{remark}
\theoremstyle{question}

\newcommand{\bu}{\mathbf{u}}

\newcommand{\bU}{\mathbf{U}}
\newcommand{\bF}{\mathbf{F}}
\newcommand{\bG}{\mathbf{G}}
\newcommand{\bR}{\mathbf{R}}
\newcommand{\jph}{j+\frac{1}{2}}
\newcommand{\jmh}{j-\frac{1}{2}}
\newcommand{\iph}{i+\frac{1}{2}}
\newcommand{\imh}{i-\frac{1}{2}}
\newcommand{\dx}{\Delta x}
\newcommand{\dy}{\Delta y}
\newcommand{\dt}{\Delta t}
\newcommand{\bau}{\bar{u}}
\newcommand{\R}{\mathbb{R}}
\newcommand{\hf}{\hat{f}}
\newcommand{\hg}{\hat{g}}

\newcommand{\ds}{\displaystyle}

\title{Non-Ocillatory Limited-Time Integration for Conservation Laws and Convection-Diffusion Equations}
\author{Jingcheng Lu and James D.Baeder}
\date{}

\begin{document}

\maketitle

\begin{abstract}
 In this study we consider unconditionally non-oscillatory, high order implicit time marching based on time-limiters.  The first aspect of our work is to propose the high resolution Limited-DIRK3 (L-DIRK3) scheme for conservation laws and convection-diffusion equations in the method-of-lines framework. The scheme can be used in conjunction with an arbitrary high order spatial discretization scheme such as 5th order WENO scheme. It can be shown that the strongly S-stable DIRK3 scheme is not SSP and may introduce strong oscillations under large time step. To overcome the oscillatory nature of DIRK3, the key idea of L-DIRK3 scheme is to apply local time-limiters (K.Duraisamy, J.D.Baeder, J-G Liu), with which the order of accuracy in time is locally dropped to first order in the regions where the evolution of solution is not smooth. In this way, the monotonicity condition is locally satisfied, while a high order of accuracy is still maintained in most of the solution domain. For convenience of applications to systems of equations, we propose a new and simple construction of time-limiters which allows flexible choice of reference quantity with minimal computation cost. Another key aspect of our work is to extend the application of time-limiter schemes to multidimensional problems and convection-diffusion equations. Numerical experiments for scalar/systems of equations in one- and two-dimensions confirm the high resolution and the improved stability of L-DIRK3 under large time steps. Moreover, the results indicate the potential of time-limiter schemes to serve as a generic and convenient methodology to improve the stability of arbitrary DIRK methods.     
\end{abstract}

\section{Introduction}
We consider high resolution, non-oscillatory implicit time integration schemes for hyperbolic conservation laws and convection-diffusion problems. The most common approach of finding numerical solutions to time-dependent partial differential equations is the \emph{method-of-lines}. In this framework, we perform proper spatial discretization over the suitable domain and then integrate the resulting system of ordinary differential equations (ODEs) using standard time integration schemes. In many practical applications this may result in a stiff ODE system, implicit time integration may be preferred in order to enlarge the allowable time step.

Most of the existing studies on non-oscillatory numerical schemes focus on spatial discretizations. The early stage framework was the class of monotone schemes, which, however, are only first order accurate \cite{monotone}. Based on less restrictive stability conditions, such as Total Variation Diminishing (TVD) and essentially non-oscillatory, many high order non-oscillatory spatial reconstructions have been proposed, such as the TVD schemes \cite{TVD} and the ENO/WENO schemes \cite{WENO}. However, ensuring non-oscillatory solutions of these high order schemes may require a severe restriction on the time step. Particularly, for high order implicit schemes this restriction can be much more restrictive than their linear stability limits.

To enable the use of larger time step for high order implicit time marching, K.Duraisamy, J.D.Baeder and J-G.Liu \cite{timelimiter} proposed the use of \emph{time-limiters}, with which the strong oscillations of trapezoidal scheme and DIRK2 scheme were effectively reduced under large time steps. To better preserve the high resolution of high order spatial reconstructions, in the present work we will extend the idea of limited-time integration to construct the new Limited-DIRK3 (L-DIRK3). The key ingredient is to locally employ implicit Euler near discontinuities. A series of numerical examples are presented to verify the improved stability and the high resolution of L-DIRK3.

\section{Motivation of Time-Limiter Schemes}
Consider initial value problem for conservation law
\begin{equation}\label{eq:CLS}
    \frac{\partial}{\partial t}u(x,t)+\frac{\partial}{\partial x}f(u(x,t)) = 0, \quad u(x,0) = u_{0}(x), \quad x\in \Omega,\quad t\in\R_{+},
\end{equation}
where $u$ is a scalar/vector of conservative variables, $f$ is a convective flux. To numerically solve (\ref{eq:CLS}), we consider the semidiscrete finite difference/finite volume schemes: at first we perform proper spatial discretization, which results in a set of ordinary differential equations 
\begin{equation}\label{eq:semidiscrete scheme}
    \frac{d}{dt}u_{j}(t)+\frac{\hf_{\jph}-\hf_{\jmh}}{\dx} = 0,\quad j = 1,2,\cdots,N.
\end{equation}
Here $N$ is the number of mesh points, $u_{j}(t)$ is the approximate solution to the point value $\displaystyle u(x_{j},t)$ /the cell average
$\displaystyle \bau(x_{j},t) : = \frac{1}{\dx}\int_{x_{\jmh}}^{x_{\jph}}u(x,t) dx$, $\displaystyle \hf_{\jph}$ is the numerical flux at cell interface $x_{\jph}$. Then we integrate (\ref{eq:semidiscrete scheme}) with standard time integration schemes, e.g. Runge Kutta methods, linear multistep methods. While explicit time integration is convenient to implement, in certain applications, e.g. steay-state computations and convection-diffusion problems, one may prefer implicit time integration in order to apply a large time step and improve the efficiency.

For nonlinear conservation laws, the solutions may develop discontinuities even if the initial data is smooth. To prevent numerical oscillations near discontinuities, many high order non-oscillatory schemes have been proposed in the last few decades and successfully applied in hyperbolic problems, including the UNO schemes \cite{UNO}, the MP schemes \cite{MP}, and the ENO/WENO type schemes \cite{WENO}, etc. However, these high order schemes are non-oscillatory only under restrictive time steps. Gottlieb, Shu and Tadmor \cite{SSP} have shown that  when the order of accuracy in time is higher than one, a time integration method, even if implicit, is at most conditionally \emph{strong stability preserving} (SSP). The time step restriction deteriorates the purpose of using of implicit methods in practical applications.

In order to use larger time step for high order implicit time integration  without introducing severe oscillations, K.Duraisamy, J. Baeder and J-G. Liu \cite{timelimiter} proposed the time-limiter schemes. These schemes are obtained from local convex combinations of a first order unconditionally SSP method and a higher order but oscillatory method. The idea is that in the regions of large solution gradients the scheme is switched to the first order method, while in the smooth regions we still apply the higher order time integration. Locally defined time-limiters are introduced to detect the occurrence of high gradients and determine the switch between different methods. Such idea was applied to construct the Limited Trapezoidal (L-Trap) and the Limited 2-Stage Diagonally Implicit Runge Kutta (L-DIRK2) methods \cite{timelimiter}. Numerical results showed that applying time-limiters effectively reduces the numerical oscillations when the time step is beyond the SSP limits of the original linear time integration schemes, and the second order accuracy is still maintained. 

In modern applications such as turbulence simulations, very high order spatial reconstructions, such as the 5th order WENO scheme, may be applied. Third order (or even higher order) implicit time integration would better preserve the high order of accuracy. An appealing option is the DIRK3 scheme, which, however, can be highly oscillatory under large CFL numbers. In the following discussions we will apply the idea of time limiting to improve the performance DIRK3.

\section{Brief Review of the DIRK3 Scheme}\label{sec:DIRK3}
In this section, we briefly review the properties of the DIRK3 scheme and explain the motivation of applying time-limiters. For convenience of notation, we represent the set of ODEs \eqref{eq:semidiscrete scheme} as a time-depedent system
\begin{equation}
    \bu_{t} = L(\bu), \quad \bu(t) = [u_{1}(t),u_{2}(t),\cdots,u_{N}(t)]^{\top}.
\end{equation}
The right hand side $L(\bu)$ represents the spatial discretization of the numerical scheme. The DIRK3 time integration scheme can be represented by the Butcher array (R. Alexander \cite{DIRK})   
\begin{equation}\label{eq:DIRK3}
\begin{tabular}{c|ccc}
       $\ds \alpha$  &  $\ds \alpha$ & 0& 0\\
       $\ds \tau_{2}$& $\ds \tau_{2}-\alpha$ & $\ds \alpha$ &0 \\
       1 & $\ds b_{1}$ &$\ds b_{2}$ & $\ds \alpha$\\
       \hline
        & $\ds b_{1}$ & $\ds b_{2}$ & $\ds \alpha$
\end{tabular}\quad, 
\end{equation}
where $\alpha \approx 0.435866521508459$ is the root of $x^{3}-3x^{2}+\frac{3}{2}x-\frac{1}{6} = 0$ lying in $(\frac{1}{6},\frac{1}{2})$,
\begin{equation*}
    \begin{split}
        & \tau_{2} = \frac{1+\alpha}{2},\\
        & b_{1} = -\frac{6\alpha^{2}-16\alpha+1}{4},\\
        & b_{3} = \frac{6\alpha^{2}-20\alpha+5}{4}.
    \end{split}
\end{equation*}
It was shown  \cite{DIRK} that \eqref{eq:DIRK3} is the unique \emph{strongly S-stable} DIRK formula of order three in three stages. The S-stability of RK methods is defined as follows.\\

\noindent \textbf{S-stability (\cite{S-stable}).} A RK method is \emph{S-stable} if for any bounded $g: [0,T]\mapsto \R$ having a bounded derivative, and any positive constant $\lambda_{0}$, there is a positive constant $h_{0}$ such that the numerical solution $\{y_{n}\}$, computed with time step $h$, to the equation 
\[
y' = g'(t)+\lambda(y-g(t))
\]
satisfies
\[
|\frac{y_{n+1}-g(t_{n+1})}{y_{n}-g(t_{n})}| <1
\]
provided $y_{n}\neq g(t_{n})$ for all $0<h<h_{0}$ all complex $\lambda$ with $Re(-\lambda)\geq\lambda_{0}$.

A RK method is \emph{strongly S-stable} if
\[
\frac{y_{n+1}-g(t_{n+1})}{y_{n}-g(t_{n})} \rightarrow 0
\]
as $Re(-\lambda)\rightarrow \infty$ for all $h>0$ such that $[t_{n},t_{n+1}] \subset [0,T]$. 

We notice that S-stability implies A-stability, which can be recovered by taking $g\equiv 0$. Being S-stable guarantees the stability of a method when applied to  problems of stiff equations. However, this does \emph{not} ensure the solutions to be non-oscillatory. In fact, the oscillatory nature of DIRK3 can be expected. We look at the first two stages of DIRK3 

stage 1:
\[
\bu^{(1)} = \bu_{n} + \alpha\dt L(\bu^{(1)}). 
\]

stage 2:
\begin{equation*}
\begin{split}
\bu^{(2)} &= \bu_{n}+(\tau_{2}-\alpha)\dt L(\bu^{(1)})+\alpha\dt L(\bu^{(2)})\\
& = \bu^{(1)}+(\tau_{2}-2\alpha)\dt L(\bu^{(1)})+\alpha\dt L(\bu^{(2)}).
\end{split}
\end{equation*}
Here the superscript of $\bu$ represents the number of stage. We notice that the second stage includes a backward time stepping of $\bu^{(1)}$ (notice that $\tau_{2}-2\alpha <0$), which is unconditionally \emph{not} TVD and may lead to oscillations in the solution. Such instability motivates the introduction of limiting mechanisms in the regions where the solution has large gradients.

\section{The Limited-DIRK3 Scheme}\label{sec:L-DIRK3}

\subsection{One-dimensional constructions}
We construct the L-DIRK3 scheme for conservation laws and convection - diffusion equations. Denote $\displaystyle \tau : = \frac{\dt}{\dx}$. The L-DIRK3 scheme for conservation equation \eqref{eq:CLS} is given by
\begin{equation}\label{eq:L-DIRK3}
    \begin{split}
        u^{(1)}_{j} & = u^{n}_{j}-\tau\alpha (\hf^{(1)}_{\jph}-\hf^{(1)}_{\jmh})\\
        u^{(2)}_{j} & = u^{n}_{j}-\tau[(a_{21,\jph}\hf^{(1)}_{\jph}-a_{21,\jmh}\hf^{(1)}_{\jmh})+(a_{22,\jph}\hf^{(2)}_{\jph}-a_{22,\jmh}\hf^{(2)}_{\jmh})]\\
        u^{n+1}_{j} &= u^{n}_{j}-\tau[(a_{31,\jph}\hf^{(1)}_{\jph}-a_{31,\jmh}\hf^{(1)}_{\jmh})+(a_{32,\jph}\hf^{(2)}_{\jph}-a_{32,\jmh}\hf^{(2)}_{\jmh})\\
        \qquad &+(a_{33,\jph}\hf^{n+1}_{\jph}-a_{33,\jmh}\hf^{n+1}_{\jmh})]
    \end{split},
\end{equation}
where 
\begin{equation*}
\begin{array}{c}
\begin{split}
    &a_{21,\jph} = \alpha+\theta^{(1)}_{\jph} (\tau_{2}-2\alpha) \\
    & a_{22,\jph} = \frac{1-\alpha}{2}+\theta^{(1)}_{\jph}(\frac{3\alpha-1}{2})\\
    & a_{31,\jph} = \alpha+\theta^{(2)}_{\jph}(b_{1}-\alpha)\\
    & a_{32,\jph} = \frac{1-\alpha}{2}+\theta^{(2)}_{\jph}(b_{2}-\frac{1-\alpha}{2})\\
    & a_{33,\jph} = \frac{1-\alpha}{2}+\theta^{(2)}_{\jph}(\frac{3\alpha-1}{2})
\end{split}
\end{array}, \quad
\theta^{(k)}_{\jph} = \frac{\theta^{(k)}_{j}+\theta^{(k)}_{j+1}}{2},\quad \theta^{(k)}_{j}\in [0,1].
\end{equation*}
The values of $\alpha$, $\tau_{2}$, $b_{1}$, $b_{2}$ are as defined in section \ref{sec:DIRK3}. The time-limiter, $\theta^{(k)}_{j}$, measures the local smoothness of solution at stage $k$. It is clear that the scheme \eqref{eq:L-DIRK3} is conservative and consistent.

To illustrate the idea of construction, we first keep $\theta^{(k)}_{j} \equiv \theta$ as constant at all points in the domain and at all stages. The resulting scheme reads 
\begin{equation}\label{eq:uniform theta}
\begin{tabular}{c|ccc}
       $\ds \alpha$  &  $\ds \alpha$ & 0& 0\\
       $\ds \tau_{2}$& $\ds \alpha+\theta (\tau_{2}-2\alpha)$ & $\ds \frac{1-\alpha}{2}+\theta(\frac{3\alpha-1}{2})$ &0 \\
       1 & $\ds \alpha+\theta(b_{1}-\alpha)$ & $\ds \frac{1-\alpha}{2}+\theta(b_{2}-\frac{1-\alpha}{2})$ & $\ds \frac{1-\alpha}{2}+\theta(\frac{3\alpha-1}{2})$\\
       \hline
        & $\ds\alpha+\theta(b_{1}-\alpha)$ & $\ds \frac{1-\alpha}{2}+\theta(b_{2}-\frac{1-\alpha}{2})$ & $\ds \frac{1-\alpha}{2}+\theta(\frac{3\alpha-1}{2})$
\end{tabular}.
\end{equation}
When $\theta = 1$ we recover the DIRK3 scheme \eqref{eq:DIRK3}, when $\theta = 0$ we obtain the unconditionally SSP successive implicit Euler steps (IE-IE-IE), 
\begin{equation}\label{eq:IE}
\begin{tabular}{c|ccc}
       $\ds \alpha$  &  $\ds \alpha$ & 0& 0\\
       $\ds \tau_{2}$& $\ds \alpha$ & $\ds \frac{1-\alpha}{2}$ &0 \\
       1 & $\ds \alpha$ & $\ds \frac{1-\alpha}{2}$ & $\ds \frac{1-\alpha}{2}$\\
       \hline
        & $\ds \alpha$ & $\ds \frac{1-\alpha}{2}$ & $\ds \frac{1-\alpha}{2}$
\end{tabular},
\end{equation}
or equivalently
\begin{equation*}
\begin{split}
    \bu^{(1)} &= \bu_{n}+\alpha\dt L(\bu^{(1)}),\\ 
    \bu^{(2)} &= \bu^{(1)}+(\tau_{2}-\alpha)\dt L(\bu^{(2)})\\
    \bu_{n+1} &= \bu^{(3)} = \bu^{(2)}+(1-\tau_{2})\dt L(\bu^{(3)})
\end{split}.
\end{equation*}
Following the idea of \cite{timelimiter}, it is expected that the first order method \eqref{eq:IE} is applied locally near discontinuities and extrema, whereas in the smooth regions we still apply the third order accurate DIRK3. Hence, the local $\theta^{(k)}_{j}$ should be applied.    

Here we propose a new and convenient construction for $\theta^{(k)}_{j}$. For scalar problems, we define
\begin{equation}\label{eq:limiter}
\begin{split}
    \theta^{(k)}_{j} &= minmod(r^{(k)}_{j},1),\\
    r^{(k)}_{j} & = \frac{u^{(k+1)}_{j+1}-u^{(k+1)}_{j-1}}{u^{(k)}_{j+1}-u^{(k)}_{j-1}}, \quad k = 1,2, \quad j = 1,2,\cdots N.
\end{split}
\end{equation}
The indicator function $\ds r^{(k)}_{j} \approx \frac{\partial_{x}u^{(k+1)}_{j}}{\partial_{x}u^{(k)}_{j}}$ is used to detect the change in the solution monotonicity at point $x_{j}$ at stage $k$. In the smooth and monotone regions we expect $r^{(k)}_{j} \approx 1$, then we would have $\theta^{(k)}_{j}\approx 1$ and  recover the DIRK3 scheme. If the solution changes the monotonicity at point $x_{j}$ when proceeds from stage $k$ to stage $k+1$ (i.e. $\partial_x u^{(k)}_j \cdot \partial_x u^{(k+1)}_j\leq 0$), which usually implies the occurrence of discontinuity or extrema, $r^{(k)}_{j}$ would be closed to 0 or become negative, then we would apply $\theta^{(k)}_{j}\approx 0$ and recover successive implicit Euler steps \eqref{eq:IE}.  

For systems of conservation laws, the concepts of monotonicity cannot be rigorously defined, the way of defining $r^{(k)}_{j}$ can be arbitrary. For exmaple, for Euler's equations we can use a density-based limiter
\[
r^{(k)}_{j} = \frac{\rho^{(k+1)}_{j+1}-\rho^{(k+1)}_{j-1}}{\rho^{(k)}_{j+1}-\rho^{(k)}_{j-1}},
\]
or a pressure-based limiter
\[
r^{(k)}_{j} = \frac{p^{(k+1)}_{j+1}-p^{(k+1)}_{j-1}}{p^{(k)}_{j+1}-p^{(k)}_{j-1}}.
\]

We point out that the L-DIRK3 scheme can be conveniently extended to convection - diffusion equations taking the form
\begin{equation}\label{eq:convec-diffu}
    \frac{\partial u}{\partial t}+\frac{\partial}{\partial x}f(u) = \frac{\partial}{\partial x}Q(u,u_x).
\end{equation}
Here the dissipation flux, $Q$, satisfies the parabolicity condition $\partial_s Q(u,s)\geq 0$, $\forall u, s$. Indeed, equation \eqref{eq:convec-diffu} admits a conservative form
\begin{equation*}
    \frac{\partial u}{\partial t}+\frac{\partial \tilde{f}}{\partial x} = 0, \quad \tilde{f}:= f(u)-Q(u,u_x).
\end{equation*}
The scheme \eqref{eq:L-DIRK3} can be naturally applied with respect to the modified flux $\tilde{f}$.

\subsection{Multidimensional extensions}
We discuss the formulation of the L-DIRK3 scheme in multi-dimensions. Without the loss of generality, we consider the two-dimensional conservation laws
\begin{equation}\label{eq:CLS 2D}
	\frac{\partial}{\partial t} u(x,y,t)+\frac{\partial}{\partial x}f(u(x,y,t))+\frac{\partial}{\partial y}g(u(x,y,t)) = 0.
\end{equation}

The conservative semidiscrete \emph{finite difference} schemes for \eqref{eq:CLS 2D} can be represented in the form
\begin{equation}\label{eq:semidiscrete 2D}
	\frac{d}{dt} u_{i,j}(t) +\frac{\hf_{\iph,j}(t)-\hf_{\imh,j}(t)}{\dx}+\frac{\hg_{i,\jph}(t)-\hg_{i,\jmh}(t)}{\dy} = 0,
\end{equation}
where $u_{i,j}(t)$ is the approximate solution to $u(x_{i},y_{j},t)$, $\hf_{\iph,j}(t)$ and $\hg_{i,\jph}(t)$ are the $x-$ and $y-$numerical fluxes such that
\begin{equation*}
	\begin{split}
		\frac{\hf_{\iph,j}(t)-\hf_{\imh,j}(t)}{\dx} &= \frac{\partial f}{\partial x}(u(x_{i},y_{j},t))+O(\dx^{m}),\\
		\frac{\hg_{i,\jph}(t)-\hg_{i,\jmh}(t)}{\dy} &= \frac{\partial g}{\partial y}(u(x_{i},y_{j},t))+O(\dy^{m}),
	\end{split}
\end{equation*}
$m$ is the desired order of accuracy in space. 

The construction of L-DIRK3 in one dimension can be extended in a dimension-by-dimension manner with slight modifications on the time limiter. To be specific, denote $\tau_{x}=\frac{\dt}{\dx}$ and $\tau_{y}=\frac{\dt}{\dy}$, the \emph{finite difference} formulation of L-DIRK3 reads
\begin{equation*}
	\begin{split}
		u^{(1)}_{i,j} & = u^{n}_{i,j}-\tau_{x}\alpha (\hf^{(1)}_{\iph,j}-\hf^{(1)}_{\imh,j})-\tau_{y}\alpha (\hg^{(1)}_{i,\jph}-\hg^{(1)}_{i,\jmh})\\
		%----------------------------------------------------------------------------
		u^{(2)}_{i,j} & = u^{n}_{i,j}-\tau_{x}(a^{x}_{21,\iph,j}\hf^{(1)}_{\iph,j}-a^{x}_{21,\imh,j}\hf^{(1)}_{\imh,j})-\tau_{y}(a^{y}_{21,i,\jph}\hg^{(1)}_{i,\jph}-a^{y}_{21,i,\jmh}\hg^{(1)}_{i,\jmh})\\
		& \quad -\tau_{x}(a^{x}_{22,\iph,j}\hf^{(2)}_{\iph,j}-a^{x}_{22,\imh,j}\hf^{(2)}_{\imh,j})-\tau_{y}(a^{y}_{22,i,\jph}\hg^{(2)}_{i,\jph}-a^{y}_{22,i,\jmh}\hg^{(2)}_{i,\jmh}),\\
		%------------------------------------------------------------------------------
		u^{n+1}_{i,j} & = u^{n}_{i,j}-\tau_{x}(a^{x}_{31,\iph,j}\hf^{(1)}_{\iph,j}-a^{x}_{31,\imh,j}\hf^{(1)}_{\imh,j})-\tau_{y}(a^{y}_{31,i,\jph}\hg^{(1)}_{i,\jph}-a^{y}_{31,i,\jmh}\hg^{(1)}_{i,\jmh})\\
		& \quad -\tau_{x}(a^{x}_{32,\iph,j}\hf^{(2)}_{\iph,j}-a^{x}_{32,\imh,j}\hf^{(2)}_{\imh,j})-\tau_{y}(a^{y}_{32,i,\jph}\hg^{(2)}_{i,\jph}-a^{y}_{32,i,\jmh}\hg^{(2)}_{i,\jmh})\\
		& \quad -\tau_{x}(a^{x}_{33,\iph,j}\hf^{n+1}_{\iph,j}-a^{x}_{33,\imh,j}\hf^{n+1}_{\imh,j})-\tau_{y}(a^{y}_{33,i,\jph}\hg^{n+1}_{i,\jph}-a^{y}_{33,i,\jmh}\hg^{n+1}_{i,\jmh}),
	\end{split}    
\end{equation*}
where
\begin{equation*}
	\begin{array}{c}
		\begin{split}
			&a^{x}_{21,\iph,j} = \alpha+\theta^{(1)}_{\iph,j} (\tau_{2}-2\alpha),  \\
			& a^{x}_{22,\iph,j} = \frac{1-\alpha}{2}+\theta^{(1)}_{\iph,j}(\frac{3\alpha-1}{2}), \\
			& a^{x}_{31,\iph,j} = \alpha+\theta^{(2),x}_{\iph,j}(b_{1}-\alpha),\\
			& a^{x}_{32,\iph,j} = \frac{1-\alpha}{2}+\theta^{(2)}_{\iph,j}(b_{2}-\frac{1-\alpha}{2}),\\
			& a^{x}_{33,\iph,j} = \frac{1-\alpha}{2}+\theta^{(2)}_{\iph,j}(\frac{3\alpha-1}{2}), 
		\end{split}
	\end{array} \hspace{0.5em} 
	\begin{array}{c}
		\begin{split}
			&a^{y}_{21,i,\jph} = \alpha+\theta^{(1)}_{i,\jph} (\tau_{2}-2\alpha),\\
			&a^{y}_{22,i,\jph} = \frac{1-\alpha}{2}+\theta^{(1)}_{i,\jph}(\frac{3\alpha-1}{2}),\\
			& a^{y}_{31,i,\jph} = \alpha+\theta^{(2),y}_{i,\jph}(b_{1}-\alpha),\\
			&a^{y}_{32,i,\jph} = \frac{1-\alpha}{2}+\theta^{(2)}_{i,\jph}(b_{2}-\frac{1-\alpha}{2}),\\
			&a^{y}_{33,i,\jph} = \frac{1-\alpha}{2}+\theta^{(2)}_{i,\jph}(\frac{3\alpha-1}{2}),
		\end{split} 
	\end{array}
\end{equation*}

\begin{equation*}
	\theta^{(k)}_{\iph.j} = \frac{\theta^{(k)}_{i,j}+\theta^{(k)}_{i+1,j}}{2}, \quad \theta^{(k)}_{i.\jph} = \frac{\theta^{(k)}_{i,j}+\theta^{(k)}_{i,j+1}}{2}.
\end{equation*}
The time limiter $\theta^{(k)}_{i,j}$ is defined as
\begin{equation*}
	\theta^{(k)}_{i,j} = minmod(r^{(k),1}_{i,j},r^{(k),2}_{i,j},r^{(k),3}_{i,j},r^{(k),4}_{i,j},1),
\end{equation*}
where $r^{(k),1}_{i,j}$, $r^{(k),2}_{i,j}$, $r^{(k),3}_{i,j}$, $r^{(k),4}_{i,j}$ are the monotonicity indicators at point $(x_i, y_j)$ from horizontal, vertical and diagonal directions
\begin{equation*}
	\begin{split}
		& r^{(k),1}_{i,j} = \frac{u^{(k+1)}_{i+1,j}-u^{(k+1)}_{i-1,j}}{u^{(k)}_{i+1,j}-u^{(k)}_{i-1,j}}, \quad
		r^{(k),2}_{i,j} = \frac{u^{(k+1)}_{i,j+1}-u^{(k+1)}_{i,j-1}}{u^{(k)}_{i,j+1}-u^{(k)}_{i,j-1}},\\
		& r^{(k),3}_{i,j} = \frac{u^{(k+1)}_{i+1,j+1}-u^{(k+1)}_{i-1,j-1}}{u^{(k)}_{i+1,j+1}-u^{(k)}_{i-1,j-1}}, \quad
		r^{(k),4}_{i,j} = \frac{u^{(k+1)}_{i-1,j+1}-u^{(k+1)}_{i+1,j-1}}{u^{(k)}_{i-1,j+1}-u^{(k)}_{i+1,j-1}}.
	\end{split}
\end{equation*}
The inclusion of information from diagonals allows for the detection of potential discontinuities that are not aligned with the axis. The extension to higher dimensions can be obtained with similar approach.

The L-DIRK3 scheme can also be applied in the framework of finite volume methods. Nevertheless, the computations will be rather complicated due to the complexity of high order quadratures for surface integrals, $\displaystyle \int^{y_{\jph}}_{y_{\jmh}}f(u(x_{\iph},y))dy$ and $\displaystyle \int^{x_{\iph}}_{x_{\imh}}g(u(x,y_{\jph}))dx$. 

\subsection{More on time-limiters} We illustrate the motivation of proposing the new limiter \eqref{eq:limiter}. In the previous work by Duraisamy et al \cite{timelimiter}, the time-limiter was constructed as
\begin{equation}\label{eq:old time-limiter}
	\theta_{j} = minmod(\frac{2s^{n+\frac{1}{2}}}{L_{j}(\bu^{n})+\epsilon}, \frac{2s^{n+\frac{1}{2}}}{L_{j}(\bu^{n+1})+\epsilon}, 1), 
\end{equation}
where $\epsilon$ is a small positive number (say $\epsilon = 10^{-10}$) introduced to enhance computational stability, and  
\begin{equation*}
	s^{n+\frac{1}{2}}_{j} = \frac{u^{n+1}_{j}-u^{n}_{j}}{\dt}, \quad L_{j}(\bu^{n})  = -\frac{\hf^{n}_{\jph}-\hf^{n}_{\jmh}}{\dx}. 
\end{equation*}
The construction \eqref{eq:old time-limiter} was motivated by Huynh's lemma \cite{Huynh} for parabolic interpolation, which implies that the solution $u_j(t)$ is monotone in the time section $t\in[t^n,t^{n+1}]$ if the two ratios in \eqref{eq:old time-limiter} are positive. 

In scalar problems, all the quantities involved in the construction \eqref{eq:old time-limiter} are computed as part of the solution update. The extra computation cost of constructing limiter is minimal. For systems of equations, one may take an arbitrary variable to construct the time-limiter. Ideally, it is preferred to use the quantity which is expected to encounter the most significant jumps. The optimal reference variable, however, may not be taken from conservative variables, such as pressure and velocity for Euler equations. In this case, the evaluations of time derivative function, $L_j$, may require expensive extra work in complicated problems. In comparison, our new construction \eqref{eq:limiter} only involves evaluations of variables and hence allows an arbitrary choice of reference quantity with minimal influence on efficiency.

\subsection{TVD global time limiting} 

In the previous study \cite{timelimiter}, primary linear monotonicity analysis was provided for local time-limiters as applied to the trapezoidal scheme. However, the rigorous proof of non-oscillatory property in general nonlinear problems is still open. To ensure the theoretical TVD-stability, a possible approach is to employ global time limiting (with uniform $\theta^{(k)}_j \equiv \theta$) proposed in the work \cite{SSPTR-BDF2} regarding the SSP variant of TR-BDF2. We define the global sensor function
\begin{equation*}
    \sigma = s_g(\bu^{n+1}) = \left\{\begin{array}{ll}
        1 &  \text{if $||\bu^{n+1}||_{TV} \geq ||\bu^n||_{TV}$} \\
        0 & \text{otherwise}
    \end{array}\right..
\end{equation*}
At each time step, tentative solution $\bu^*$ is obtained with the unlimited DIRK3 ($\theta = 1$). If $\sigma = s_g(\bu^*) = 0$, we accept $\bu^*$ as the new solution. If $\sigma = 0$, then the tentative solution violates the stability constraint and we would repeat the integration with implicit Euler ($\theta = 0$). In this way, the limited scheme naturally inherits the $A-$stability of DIRK3 and implicit Euler. The TVD-stability is also ensured, since implicit Euler is unconditionally TVD and DIRK3 is applied only if the scheme does not increase the total variation.

However, when solving non-smooth problems the approach of global time limiting suffers from the risk of doubling computation cost. Indeed, the unlimited DIRK3 may be too oscillatory such that the time integration is restarted frequently. Meanwhile, the use of IE-IE-IE mode may result in order reduction to first order, which was reported in \cite{SSPTR-BDF2}. Hence in practical applications, local time limiting would be recommended to ensure good accuracy and high efficiency. In fact, the numerical results in section \ref{sec:numeric} will show that applying local time-limiters can effectively enhance the stability of DIRK3 under large time step.

\section{Numerical Examples}\label{sec:numeric}
We present numerical results for several test cases. For all the computations of DIRK3 and L-DIRK3, we apply the 5th order WENO-JS finite difference scheme \cite{WENO} for space reconstruction. The number of equally spaced mesh points in the domain is represented by $N$. The implicit systems are solved with Newton-type sub-iterations \cite{Newton}. At each implicit stage, the iterations stop if the $l^2-$norm of the residual is reduced by four orders of magnitudes, or the number of iterations reaches the maximum of 30. The limiters $\theta^{(k)}_{j}$ are updated explicitly in the iterations. Verification of the accuracy of WENO solver was presented in the earlier manuscript \cite{ldirk3scitech}.

\subsection{Linear advection equation}

\begin{figure}[h!]
	\centering
	\includegraphics[scale = 0.3]{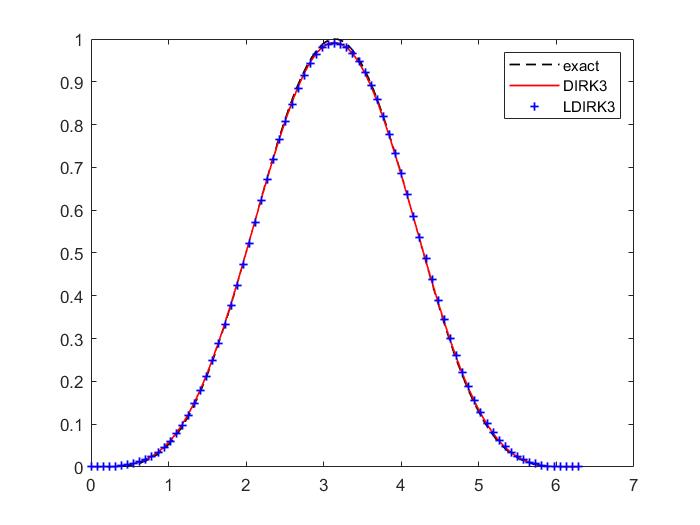}
	\caption{Linear advection, $\sin^{4}(\frac{x}{2})$ wave, DIRK3 vs L-DIRK3, $CFL = 0.5$, 5th order WENO in space, $N = 80$, periodic bc, 1 period of evolution.}
	\label{fig:sinwave}
\end{figure}

The first test case is the linear advection equation 
\begin{equation*}
	\frac{\partial u}{\partial t}+\frac{\partial u}{\partial x} = 0
\end{equation*}
with periodic boundary conditions and a smooth initial condition $\ds u_{0}(x) = \sin^{4}(\frac{x}{2})$ over the computational domain $[0, 2\pi]$. This test case is chosen to demonstrate the order of accuracy in time and the resolution at smooth extrema. Tables \ref{tab: linear err ldirk3}--\ref{tab: linear err dirk3} present the $L^{\infty}-$, $L^{1}-$ and $L^{2}-$errors of different methods under $N = 25, 50, 100, 200, 400$ at $CFL = 0.5$ after one period of evolution. The results indicate the second-order convergence rates of L-DIRK3. The uniform third order accuracy in time is not recovered since the scheme reduces to first order at the extrema. Figure \ref{fig:sinwave} shows the solutions by L-DIRK3 and DIRK3 under $N = 80$. We observe that L-DIRK3, although suffers from order reduction, resolves the extrema as sharply as the unlimited DIRK3 under a moderate mesh size.

\begin{table}[h!]
\renewcommand\arraystretch{1.25} 
\centering
\begin{tabular}{ccccccc} 
\hline
\hline
N & $L^{\infty}$ err &  $L^{1}$ err & $L^2$err & $L^{\infty}$ rate & $L^{1}$ rate & $L^{2}$ rate   \\
\hline
50& $1.58e^{-2}$ & $2.19e^{-3}$ & $3.67e^{-3}$&-- & -- & --\\
100& $4.72e^{-3}$ & $6.91e^{-4}$ & $1.05e^{-3}$& 1.74 & 1.66 & 1.81 \\
200 & $1.24e^{-3}$ & $1.87e^{-4}$ & $2.78e^{-4}$&1.93 & 1.88 & 1.91\\
400 & $3.85e^{-4}$ & $4.60e^{-5}$ & $6.89e^{-5}$& 1.68 & 2.03 & 2.01\\
800 & $1.09e^{-4}$ & $1.02e^{-5}$ &$1.66e^{-5}$ & 1.82 & 2.17 & 2.06 \\
\hline
\hline
\end{tabular}
\caption{Error norms for linear advection, L-DIRK3, CFL = 0.5, 5th order WENO in space, periodic bc, 1 period of evolution.}  
\label{tab: linear err ldirk3}  
\end{table}

\begin{table}[h!]
	\renewcommand\arraystretch{1.25} 
	\centering
	\begin{tabular}{ccccccc} 
		\hline
		\hline
		N & $L^{\infty}$ err &  $L^{1}$ err & $L^2$err & $L^{\infty}$ rate & $L^{1}$ rate & $L^{2}$ rate   \\
		\hline
		50& $6.07e^{-4}$ & $1.96e^{-4}$ & $2.37e^{-4}$&-- & -- & --\\
		100& $1.71e^{-5}$ & $9.01e^{-6}$ & $1.01e^{-5}$& 5.15 & 4.44 & 3.55 \\
		200 & $1.71e^{-6}$ & $8.16e^{-7}$ & $9.53e^{-7}$&3.32 & 3.46 & 3.40\\
		400 & $2.00e^{-7}$ & $1.02e^{-7}$ & $1.16e^{-7}$& 3.09 & 3.00 & 3.04\\
		800 & $2.47e^{-8}$ & $1.27e^{-8}$ &$1.44e^{-8}$ & 3.02 & 3.00 & 3.01\\
		\hline
		\hline
	\end{tabular}
	\caption{Error norms for linear advection, DIRK3, CFL = 0.5, 5th order WENO in space, periodic bc, 1 period of evolution.}  
	\label{tab: linear err dirk3}  
\end{table}

\begin{figure}[h!]	
	\centering
	\includegraphics[scale = 0.3]{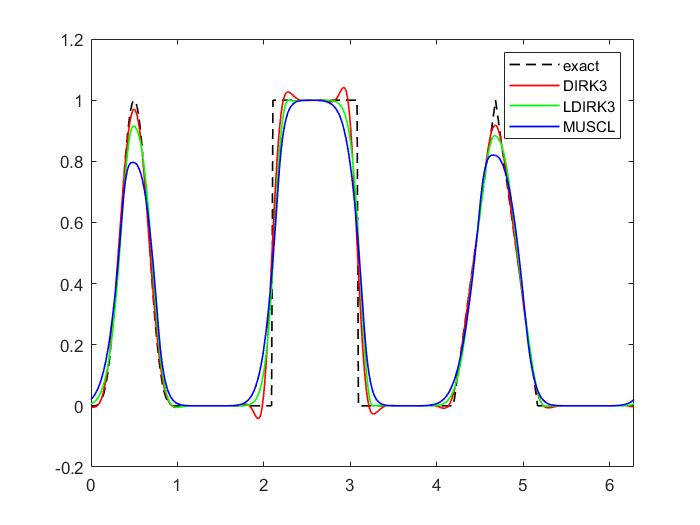}
	\caption{Linear advection, multiple waves, DIRK3 vs L-DIRK3, $CFL = 2$, 5th order WENO in space, $N = 400$, periodic bc, 1 period of evolution.}
	\label{fig:multiwave}
\end{figure}

The second test case comprises of multiple waves over the domain $[0, 2\pi]$, including a $\sin^{4}(\pi x)$ distribution, a step function and a hat function. The mesh size $N = 400$ is applied. Figure \ref{fig:multiwave} shows the solutions by DIRK3 and L-DIRK3 at $CFL = 2$ after one period of evolution. We also include the second order solution for comparison, which is computed with piecewise linear MUSCL extrapolation \cite{MUSCL} (with minmod slope limiter) in space and explicit SSPRK3 \cite{SSP} time integration at $CFL = 0.9$. It is seen that the solution by DIRK3 presents obvious overshoots and undershoots near the edges of square wave, whereas the L-DIRK3 generates non-oscillatory solution. In spite of slight clipping at the tops of sine wave and triangular wave, the solution by L-DIRK3 resolves the extremum obviously better than the second order solution.

\subsection{1D Burger's equation}
\begin{figure}[h!]
	\centering
	\includegraphics[scale = 0.3]{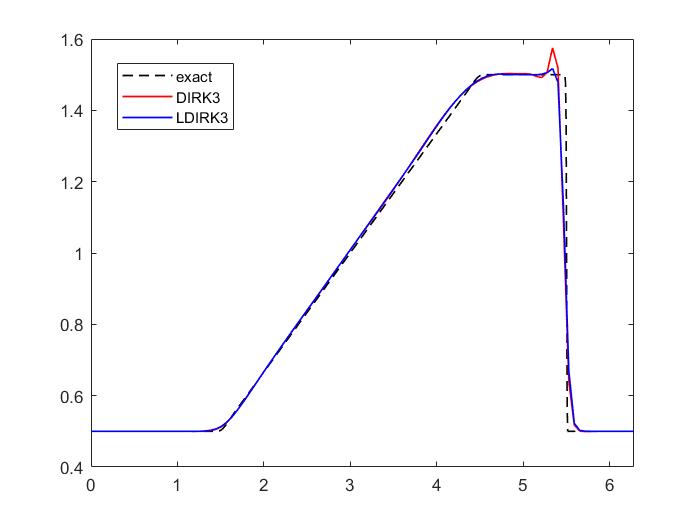}
	\caption{1D Burger's equation, expansion and comrpession wave, DIRK3 vs L-DIRK3, $CFL = 3$, 5th order WENO in space, $N = 100$, $t = 2$.}
	\label{fig:exp comp LDIRK3}
\end{figure}

We present the results for the one-dimensional Burger's equation
\begin{equation*}
	\frac{\partial u}{\partial t}+\frac{\partial}{\partial x}(\frac{u^{2}}{2}) = 0
\end{equation*} 
with periodic boundary conditions over the domain $[0, 2\pi]$. The initial data comprises of an expansion wave and a compression wave. Figure \ref{fig:exp comp LDIRK3} shows the solutions at $t = 2$ with mesh size $N = 100$ at $CFL = 3$. The solution comprises of an expansion wave and a shock wave. It is seen that the solution by DIRK3 has obvious overshoot near the shock under relatively large time step. The solution by L-DIRK3 is much less oscillatory and resolves the expansion wave and the shock well.

\subsection{1D Euler equations}

The one-dimensional Euler's equations of gas dynamics is given by
\begin{equation*}
    \frac{\partial \bU}{\partial t}+\frac{\partial \bF}{\partial x} = 0,
\end{equation*}
where $\bU$, the vector of conserved variables, and $\bF$, the flux vector, are given by
\begin{equation*}
    \bU = \left[\begin{array}{c}
         \rho  \\
         \rho u \\
         e
    \end{array}\right], \qquad
    \bF = \left[\begin{array}{c}
         \rho u \\
         \rho u^{2}+p \\
         (e+p)u
    \end{array}\right],
\end{equation*}
$\rho$, $u$, $p$ are density, velocity and pressure, $e$ is the total energy per unit volume which is given by
\[
e = \frac{p}{\gamma-1}+\frac{\rho u^{2}}{2}, \quad \gamma = 1.4.
\]

We will test the performance of L-DIRK3 in the Sod's problem, Lax's problem and Osher-Shu problems. These test cases are Riemann problems in a constant area tube. The left and the right states are indexed by $L$ and $R$. 

The interface numerical fluxes are constructed with the characteristic-wise Lax Friedrich flux vector splitting 
\begin{equation}\label{eq:char LF}
	\begin{split}
		\bF^{\pm}(\bU) &= \frac{1}{2}(\bF(\bU)\pm \bR\Lambda \bR^{-1}\bU), \quad \Lambda = diag(\max_{\bU}|\{\lambda_{i}(\bU)|\}_{i}),\\
		\hat{\bF}_{\jph} &=  \hat{\bF}^{+}_{\jph}+\hat{\bF}^{-}_{\jph}.
	\end{split}    
\end{equation}
Here $\bR$ is the matrix of right eigenvectors of flux Jacobian $\ds \frac{\partial \bF}{\partial \bU}$, $\hat{\bF}_{\jph}$ is the numerical flux at cell interface $x_{\jph}$, $\hat{\bF}^{+}_{\jph}$, $\hat{\bF}^{-}_{\jph}$ are the reconstructions of $\bF^{+}$, $\bF^{-}$ at $x_{\jph}$ from the left and the right.

\paragraph{Sod's problem} 

\begin{figure}[h!]
    \centering
    \begin{subfigure}{0.35\textwidth}
    \includegraphics[scale = 0.25]{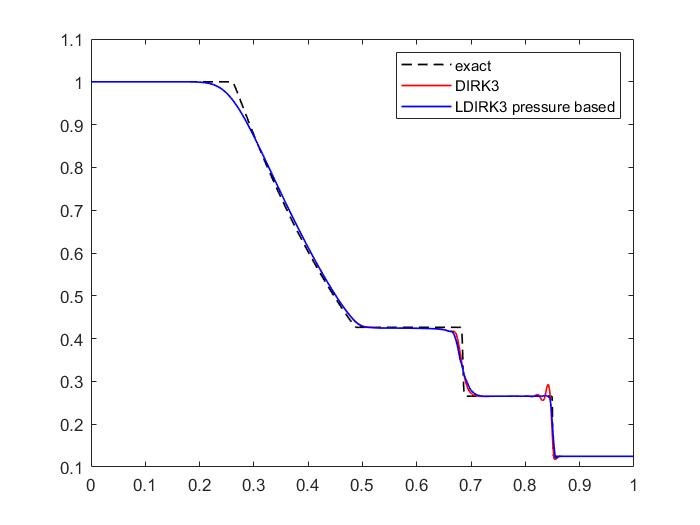}
    \subcaption{Pressure-based}
    \end{subfigure}
    \quad
    \begin{subfigure}{0.35\textwidth}
    \includegraphics[scale = 0.25]{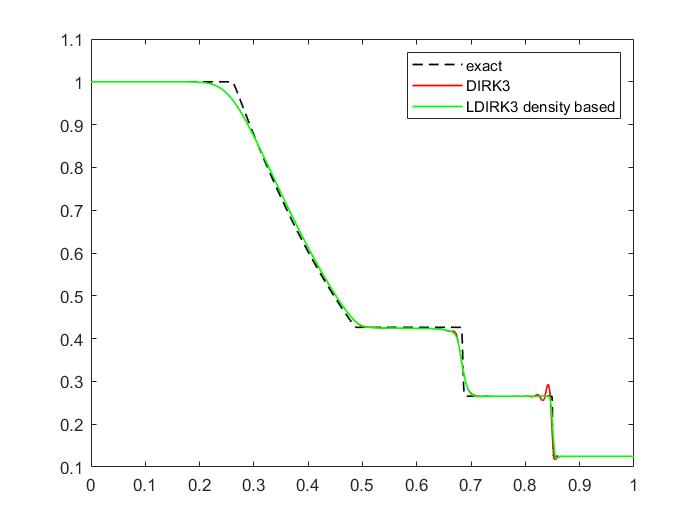}
    \subcaption{Density-based}
    \end{subfigure}
    \caption{Sod's problem (Density), $CFL = 4$, 5th order WENO in space, $N=400$, $t = 0.2$.}
    \label{fig:ST}
\end{figure}

The computational domain of Sod's problem is [0, 1] and the interface is at $x=0.5$. The initial data is given by
\begin{equation*}
\left\{\begin{array}{c}
\begin{split}
        & p_{L} = 1.0, \quad p_{R} = 0.1  \\
        & \rho_{L} = 1.0, \quad \rho_{R} = 0.125 \\
        & u_{L} = 0.0, \quad u_{R} = 0.0 
\end{split}
\end{array}\right..
\end{equation*}
The numerical solution is computed with $N = 400$ under $CFL = 4$. We present the evolution of density at $t = 0.2$ in Figure \ref{fig:ST}. We see that under high CFL number the solution of DIRK3 has obvious oscillation at the shock. By comparison, the solutions by L-DIRK3, with density-based and pressure-based limiting,  have no visible oscillation at the shock, while the resolution at the expansion wave and discontinuities remain comparable with the unlimited DIRK3. 

\paragraph{Lax's problem}

\begin{figure}[h!]
\vspace{0.2in}
    \centering
    \begin{subfigure}{0.35\textwidth}
    \includegraphics[scale = 0.25]{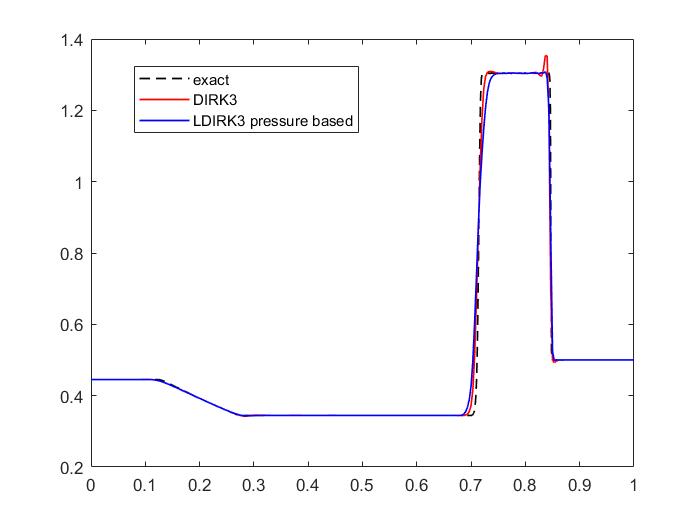}
    \subcaption{Pressure-based}
    \end{subfigure}
    \quad
    \begin{subfigure}{0.35\textwidth}
    \includegraphics[scale = 0.25]{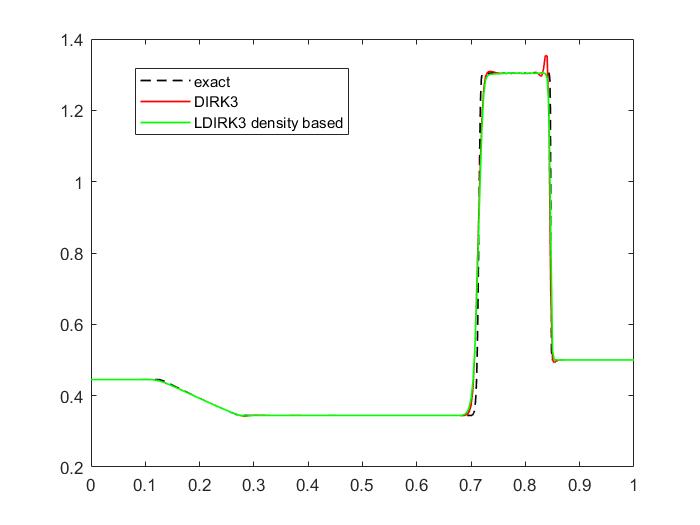}
    \subcaption{Density-based}
    \end{subfigure}
    \caption{Lax's problem (Density), $CFL = 3.5$, 5th order WENO in space, $N=300$, $t = 0.14$.}
    \label{fig:Lax's}
\end{figure}

The initial data of Lax's problem is given by
\begin{equation*}
\left\{\begin{array}{c}
\begin{split}
        & p_{L} = 3.528, \quad p_{R} = 0.571  \\
        & \rho_{L} = 0.445, \quad \rho_{R} = 0.5 \\
        & u_{L} = 0.698, \quad u_{R} = 0.0 
\end{split}
\end{array}\right..
\end{equation*}
over the domain $[0, 1]$ with interface at $x = 0.5$. We compute the solutions at $t = 0.14$ with mesh size $N = 300$ at $CFL = 3.5$. Figure \ref{fig:Lax's} shows the density evolution. It is seen that the unlimited DIRK3 generates obvious oscillation at the edge of the square wave. The solutions by L-DIRK3, both pressure-based and density- based, present no visible oscillation and resolve the expansion wave and the shocks well.

\paragraph{Osher-Shu problem}

\begin{figure}[h!]
    \centering
    \begin{subfigure}{0.35\textwidth}
    \includegraphics[scale = 0.25]{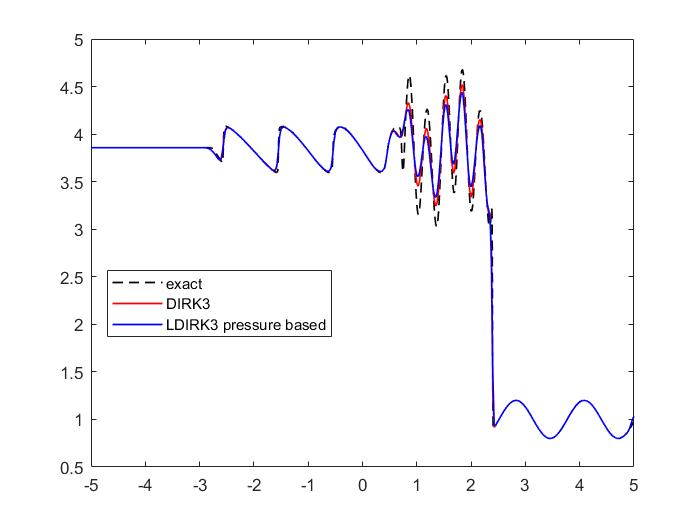}
    \subcaption{Pressure-based}
    \end{subfigure}
    \quad
    \begin{subfigure}{0.35\textwidth}
    \includegraphics[scale = 0.25]{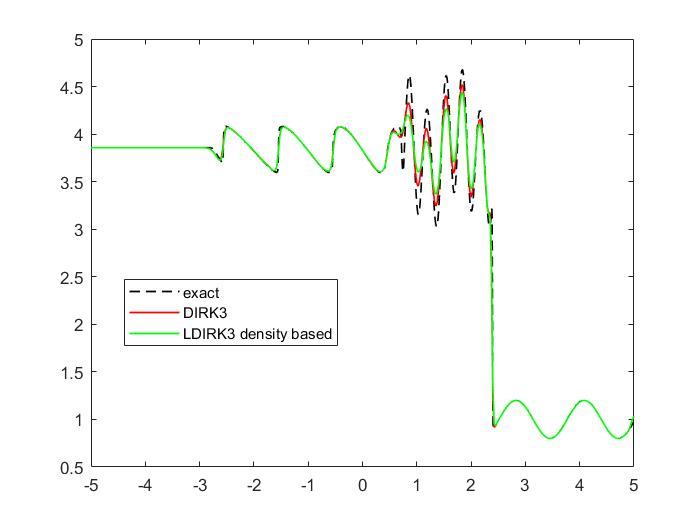}
    \subcaption{Density-based}
    \end{subfigure}
    \caption{Osher-Shu (Density), $CFL = 2$, 5th order WENO in space, $N=400$, $t = 1.8$.}
    \label{fig:SE cfl2}
\end{figure}

In Osher-Shu problem, the solution comprises of a discontinuity as well as a smooth harmonic density wave. We will examine the resolution of L-DIRK3 at extremum and its stability at the shock. The computational domain is [-5, 5] and the interface is at $x = -4$. The initial data is given by
\begin{equation*}
\left\{\begin{array}{c}
\begin{split}
        & p_{L} = 10.33333, \quad p_{R} = 1.0  \\
        & \rho_{L} = 3.857143, \quad \rho_{R} = 1+0.2\sin(5x) \\
        & u_{L} = 2.6293690, \quad u_{R} = 0.0 
\end{split}
\end{array}\right..
\end{equation*}
We present the evolution of density at $t = 1.8$. The mesh size $N = 400$ is applied. At first, we test with $CFL = 2$. The result is shown in Figure \ref{fig:SE cfl2}. We see that under relatively large CFL number, the solutions by L-DIRK3 have a little clipping at extremum but not too severe. The resolutions of the L-DIRK3 and the unlimited DIRK3 are overall comparable.

\begin{figure}[h!]
    \centering
    \begin{subfigure}{0.35\textwidth}
    \includegraphics[scale = 0.25]{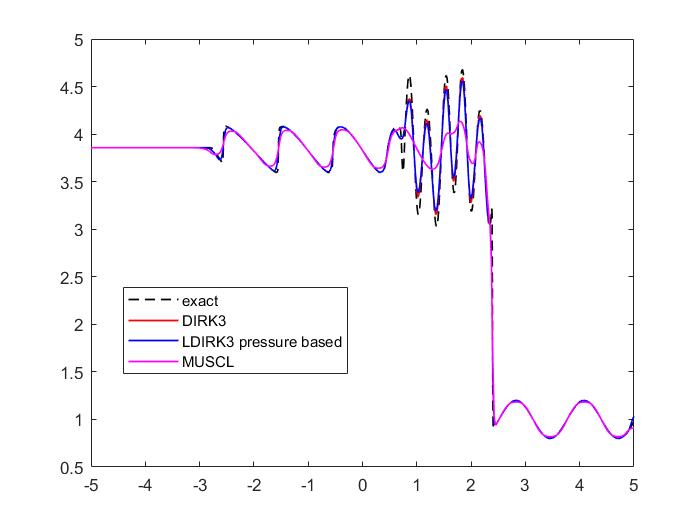}
    \subcaption{Pressure-based}
    \end{subfigure}
    \quad
    \begin{subfigure}{0.35\textwidth}
    \includegraphics[scale = 0.25]{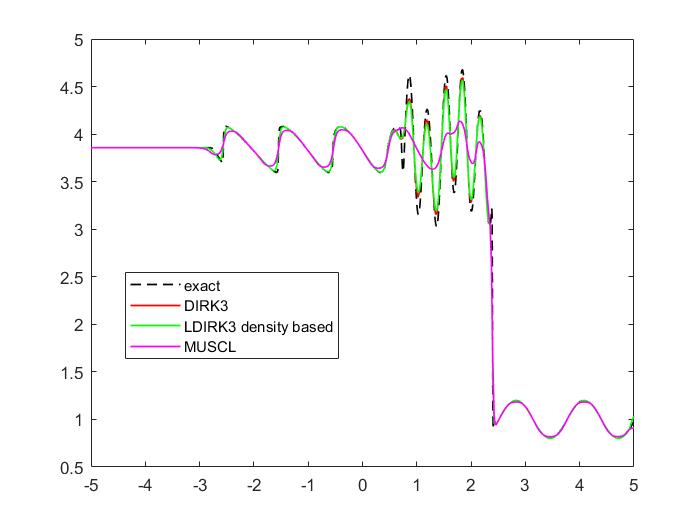}
    \subcaption{Density-based}
    \end{subfigure}
    \caption{Osher-Shu (Density), $CFL = 0.9$, 5th order WENO in space, $N=400$, $t = 1.8$.}
    \label{fig:SE cfl0.9}
\end{figure}

We also present the solutions under $CFL = 0.9$ in Figure \ref{fig:SE cfl0.9}. It is seen that when the smaller time step is applied, the solutions given by DIRK3 and L-DIRK3 become almost indistinguishable. For comparison, we include the second order solution under the same mesh size and the same CFL condition in the plot. For the second order solution, we apply the piecewise linear MUSCL extrapolation (with minmod slope limiter) for space discretization, the interface numerical fluxes are computed with the flux vector splitting \eqref{eq:char LF}, and we apply the explicit SSPRK3 for time integration. We see that the solutions by L-DIRK3, both density- and pressure-based, capture the extremum much more sharply than the typical second order solution. The high resolution of 5th order WENO scheme is well preserved despite of the time limiting.      

\begin{figure}[h!]
    \centering
    \begin{subfigure}{0.35\textwidth}
    \includegraphics[scale = 0.25]{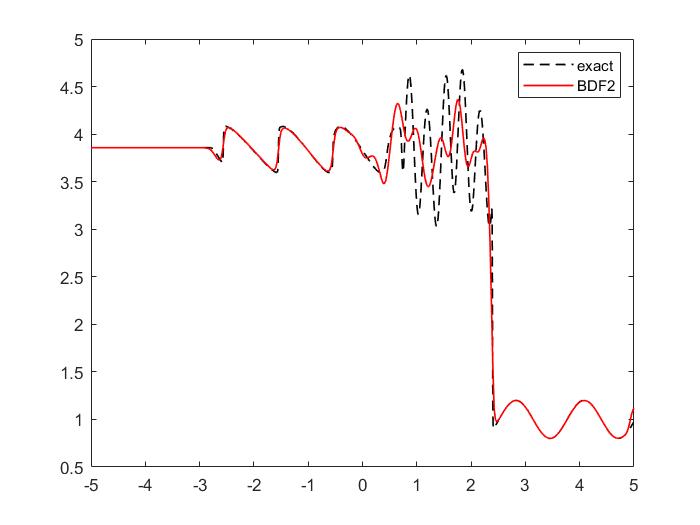}
    \subcaption{BDF2}
    \end{subfigure}
    \quad
    \begin{subfigure}{0.35\textwidth}
    \includegraphics[scale = 0.25]{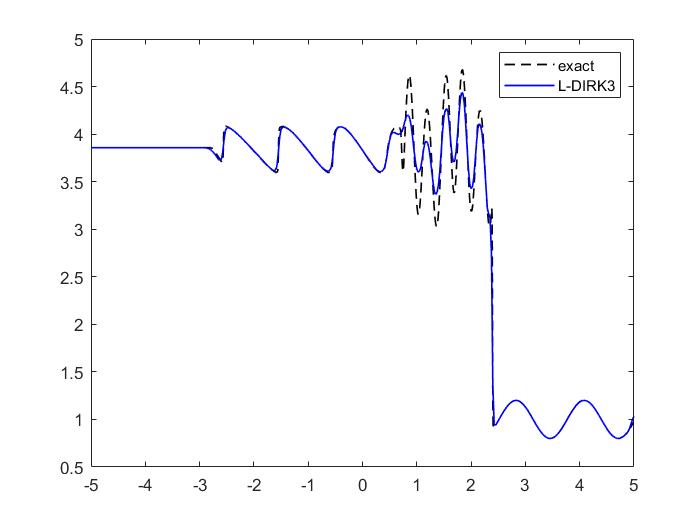}
    \subcaption{L-DIRK3}
    \end{subfigure}
    \caption{Osher-Shu (Density), BDF2 vs L-DIRK3, $CFL = 2$, 5th order WENO in space, $N=400$, $t = 1.8$.}
    \label{fig:SE BDF2 vs LDIRK3}
\end{figure}

In turbulence simulations, BDF2 is an often used second order implicit time marching method (see e.g. \cite{martin2006,jia2019}). Here we also compare the performances of L-DIRK3 and BDF2 under relatively large CFL number, as shown in Figure \ref{fig:SE BDF2 vs LDIRK3}. It is observed that BDF2 introduces notable dispersive errors while L-DIRK3 captures the wave positions much more accurately.

\subsection{Convection-diffusion equations}
In real world applications, a scenario where implicit methods may actually be preferred over explicit methods is when solving convection - diffusion problems, where the time step scales at $\dt\sim\dx^2$, in contrast to $\dt\sim\dx$ for convective problems. A large CFL number is desired in order to produce a reasonable time step. The diffusion term is often discretized implicitly to get rid of the severe time step restriction. We first consider the viscous Burger's equation 
\[
\frac{\partial u}{\partial t}+\frac{\partial}{\partial x}(\frac{u^2}{2}) = \epsilon\frac{\partial }{\partial x}(\nu(u)\frac{\partial u}{\partial x}), 
\] 
 associated with discontinuous diffusion coefficient 
\begin{equation}\label{eq:discontinuous diffusion}
\nu(u) = \left\{\begin{array}{ll}
1, \quad & |u|\geq 0.5,\\
0, \quad & \text{otherwise}.
\end{array}
\right.
\end{equation}
 The equation has hyperbolic nature when $u\in[-0.5, 0.5]$ and parabolic elsewhere. 
 
\begin{figure}[h!]
    \centering
    \begin{subfigure}{0.4\textwidth}
    \includegraphics[scale = 0.25]{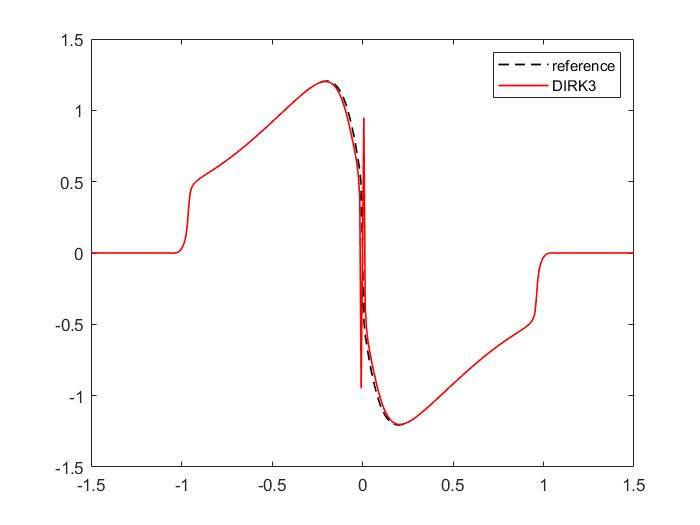}
    \subcaption{DIRK3}
    \end{subfigure}
    \quad    
    \begin{subfigure}{0.5\textwidth}
    \includegraphics[scale = 0.25]{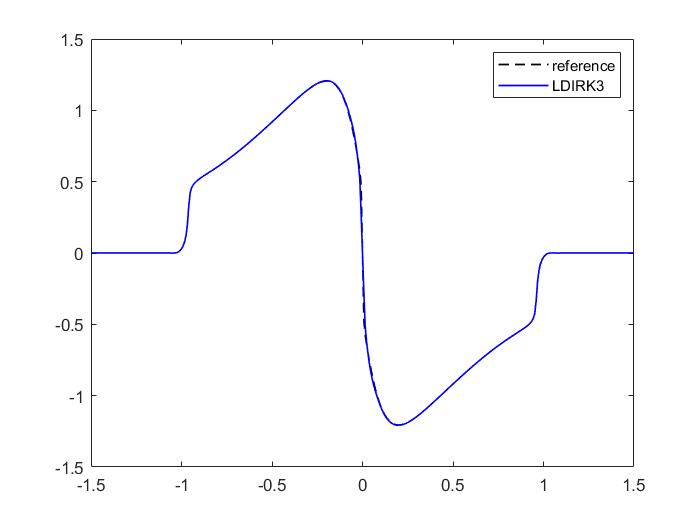}
    \subcaption{L-DIRK3}
    \end{subfigure}
    \caption{Viscous Burger's equation, DIRK3 vs L-DIRK3, $CFL = 10$, 5th order WENO in space, $N = 400$, $t = 0.5$.}
    \label{fig:ConvecDiffu1D burgers}
\end{figure}

In the following computations, the initial data are set to 2 and -2 on the intervals $[-0.9, -0.1]$ and $[0.1, 0.9]$ respectively over the domain $[-1.5, 1.5]$. the time step is chosen according to the $CFL$ condition 
\[
\dt(\frac{a}{\dx}+\max_u\frac{\epsilon \nu(u)}{\dx^2})= CFL,\quad a = \max_u |f'(u)|.
\]
To enhance the convergence of implicit iterations while maintaining the sharp transition between different flow regimes, we approximate the diffusion coefficient \eqref{eq:discontinuous diffusion} with smooth hyperbolic tangent
\[ \nu^{\epsilon_1}(u) = \frac{1}{2}(1+\tanh(\frac{|u|-0.5}{\epsilon_1})), \quad \epsilon_1 = 0.05.
\]
The convection term is discretized with fifth order WENO method and the diffusion term is discretized with the usual second order central differencing.  We take $\epsilon - 0.1$ and compute the solution at $t = 0.2$ under mesh size $N = 400$ a large $CFL$ number of 10. Figure \ref{fig:ConvecDiffu1D burgers} compares the results by DIRK3 and L-DIRK3. The reference solution is computed with explicit SSPRK3 at $CFL = 0.6$ under the same space mesh size. 
The solution by DIRK3 generates strong oscillations at the jump, whereas the L-DIRK3 scheme generates non-oscillatory solution under a much larger time step than that of the explicit method. 

\begin{figure}[h!]
    \centering
    \begin{subfigure}{0.43\textwidth}
    \includegraphics[scale = 0.25]{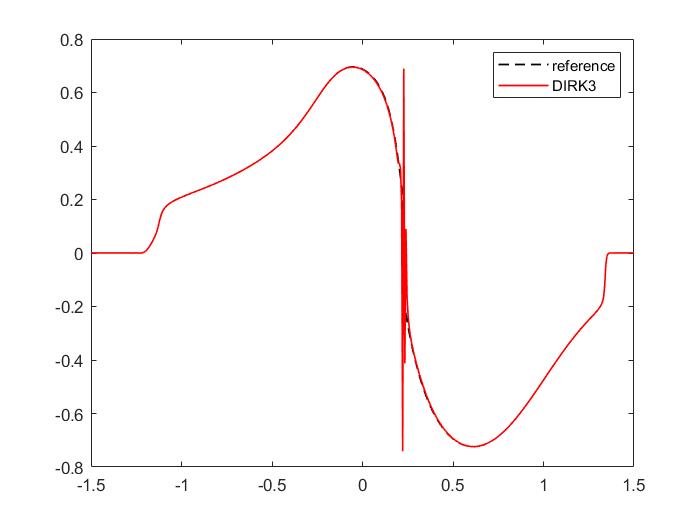}
    \subcaption{DIRK3}
    \end{subfigure}
    \quad    
    \begin{subfigure}{0.5\textwidth}
    \includegraphics[scale = 0.25]{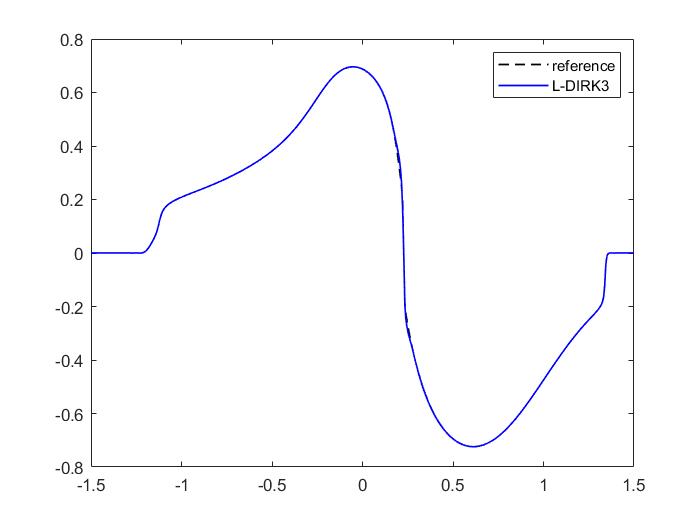} 
    \subcaption{L-DIRK3}
    \end{subfigure}
    \caption{Viscous Buckley-Leverette equation, DIRK3 vs L-DIRK3, $CFL = 10$, 5th order WENO in space, $N = 500$, $t = 0.5$.}
    \label{fig:ConvecDiffu1D B-L}
\end{figure}

Next, we consider the viscous Buckley-Leverett equation
\[
\frac{\partial u}{\partial t}+\frac{\partial}{\partial x}f(u) = \epsilon\frac{\partial }{\partial x}(\nu(u)\frac{\partial u}{\partial x}), \quad f(u) = \frac{u^2}{u^2+(1-u)^2}.
\] 
This equation is a prototype model for oil reservoir (two-phase flow) simulation. The diffusion coefficient is taken to 
\[
\nu^{\epsilon_1}(u) = \frac{1}{2}(1+\tanh(\frac{|u|-0.2}{\epsilon_1})), \quad \epsilon_1 = 0.03,
\]
which is a smooth approximation to discontinuous function $\nu(u) = 1$, $|u|\geq 0.2$. The initial data are set to 0.9 for $|x+\frac{1}{\sqrt{2}}|<0.4$  and -0.9 for $|x-\frac{1}{\sqrt{2}}|<0.4$.  We compute the solution with respect to $\epsilon = 0.1$ under $CFL = 10$, as shown in Figure \ref{fig:ConvecDiffu1D B-L}. Again, strong oscillations are observed in the solution by the unlimited DIRK3 scheme, whereas the result by L-DIRK3 presents a satisfying stability as well as a high resolution.

\vspace{0.3in}

\subsection{Two-dimensional problems}

We examine the performance of L-DIRK3 as applied to two-dimensional scalar and system of equations.

\paragraph{2D Buckley-Leverett equation}

\begin{figure}[h!]
    \centering
    \begin{subfigure}{0.45\textwidth}
    \centering
    \includegraphics[scale = 0.22]{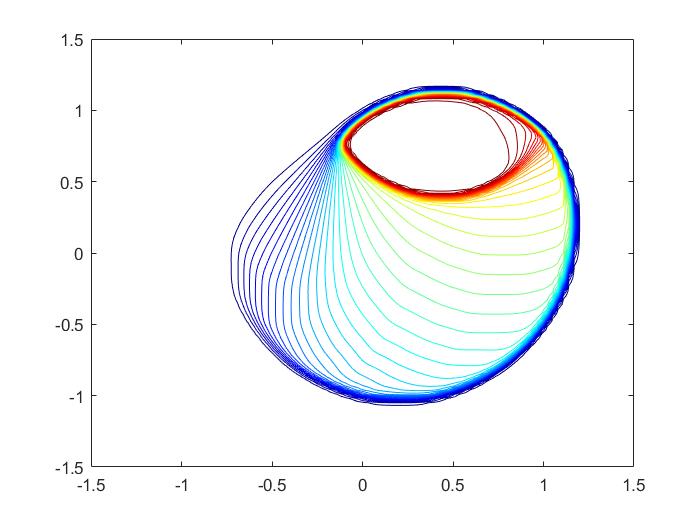}
    \subcaption{Solution contours}
    \end{subfigure}\quad
    \begin{subfigure}{0.5\textwidth}
    \centering
    \includegraphics[scale = 0.22]{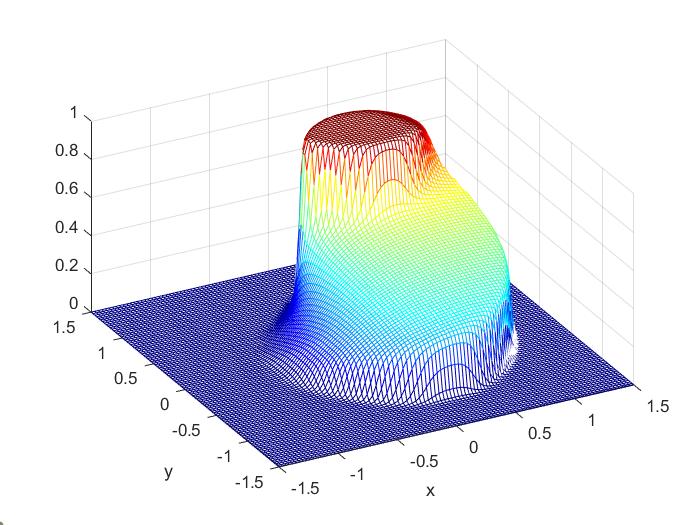}
    \subcaption{Surface plot}
    \end{subfigure}
    \caption{2D Buckley-Leverett equation, reference solution, $N_{x} = N_{y} = 100$, $t = 0.4$.}
    \label{fig:2D B-L ref}
\end{figure}

\begin{figure}[h!]
    \centering
        \begin{subfigure}{0.4\textwidth}
    \centering
    \includegraphics[scale = 0.25]{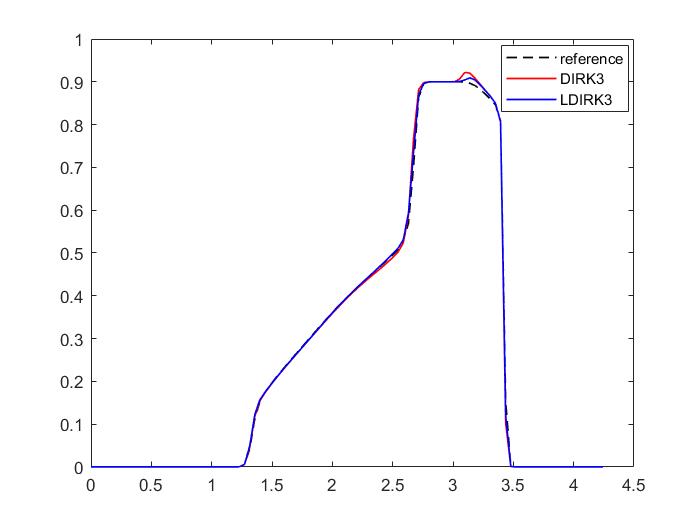}
    \subcaption{Solution contours}
    \end{subfigure}\quad
    \begin{subfigure}{0.5\textwidth}
    \centering
    \includegraphics[scale = 0.25]{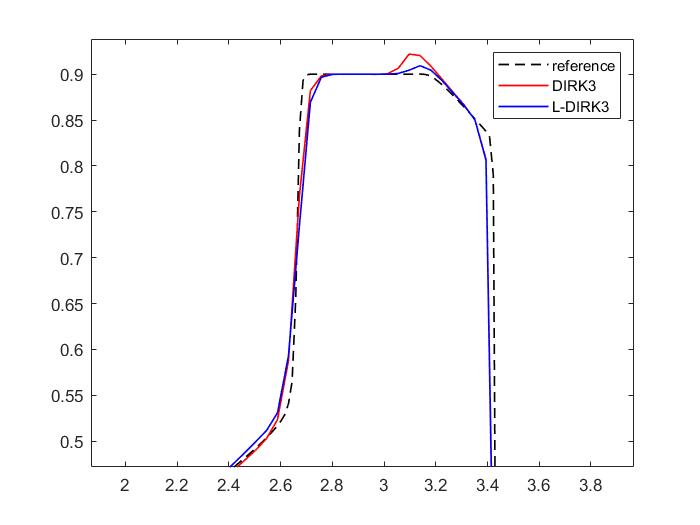}
    \subcaption{Surface plot}
    \end{subfigure}
    \caption{2D Buckley-Leverett equation, reference solution, $N_{x} = N_{y} = 100$, $t = 0.4$.}
    \label{fig:2D B-L diagonal} 
\end{figure}

Consider the two-dimensional inviscid Buckley-Leverett equation 
\begin{equation*}
	\frac{\partial u}{\partial t} + \frac{\partial}{\partial x}f(u)+\frac{\partial}{\partial y}g(u) = 0,
\end{equation*}
with the flux functions
\begin{equation*}
	\begin{split}
		f(u)  &= \frac{u^{2}}{(1-u)^{2}+u^{2}} ,\\
		g(u)  &= f(u)(1-5(1-u)^{2}),
	\end{split}    
\end{equation*}
and the initial data
\begin{equation*}
	u(x,y,0) = \left\{\begin{array}{c}
		\begin{split}
			0.9 & \hspace{1em} x^{2}+y^{2}<0.5 \\
			0 &  \hspace{1em} \text{else}
		\end{split}
	\end{array}\right..
\end{equation*}
The reference solution, which is computed with SSPRK3-WENO5 scheme at $CFL = 0.9$, is shown in Figure \ref{fig:2D B-L ref}.  Figure \ref{fig:2D B-L diagonal} compares the solutions by different methods along the diagonal $x = y$. It is seen that L-DIRK3, with the help of time limiting, has a smaller overshoot than DIRK3 near the shock. 

\paragraph{2D Euler equations}
We apply the L-DIRK3 scheme to solve the two-dimensional Euler's equations for an ideal gas
\begin{equation*}
	\frac{\partial \bU}{\partial t} + \frac{\partial \bF(\bU)}{\partial x}+\frac{\partial \bG(\bU)}{\partial y} = 0,
\end{equation*}
the conservative variables $\bU$, the $x-$flux $\bF$ and the $y-$flux $\bG$ are given by
\begin{equation*}
	\bU = \left[\begin{array}{c}
		\rho  \\
		\rho u \\
		\rho v\\
		e
	\end{array}\right], \quad
	\bF = \left[\begin{array}{c}
		\rho u  \\
		\rho u^2 +p \\
		\rho u v\\
		(e+p) u 
	\end{array}\right],\quad
	\bG = \left[\begin{array}{c}
		\rho v  \\
		\rho u v \\
		\rho v^2 + p\\
		(e+p) v 
	\end{array}\right],
\end{equation*}
where $\rho$, $p$ are density and pressure, $u$, $v$ are velocities in $x-$ and $y-$direction, the total energy $e$ is given by
\[
e = \frac{p}{\gamma-1}+\frac{\rho (u^2 +v^2)}{2}
\]
with $\gamma = 1.4$.

\begin{figure}[h!]
	\centering
	\begin{subfigure}{0.43\textwidth}
	\includegraphics[scale = 0.25]{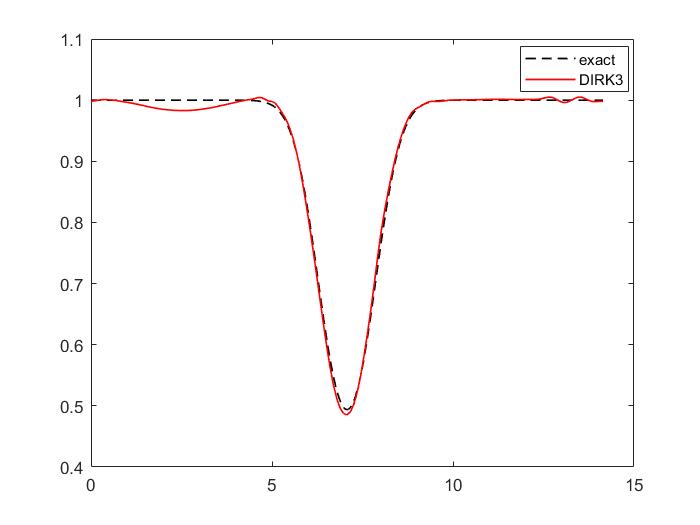}
	\subcaption{DIRK3}
	\end{subfigure}
	\quad
	\begin{subfigure}{0.5\textwidth}
	\includegraphics[scale = 0.25]{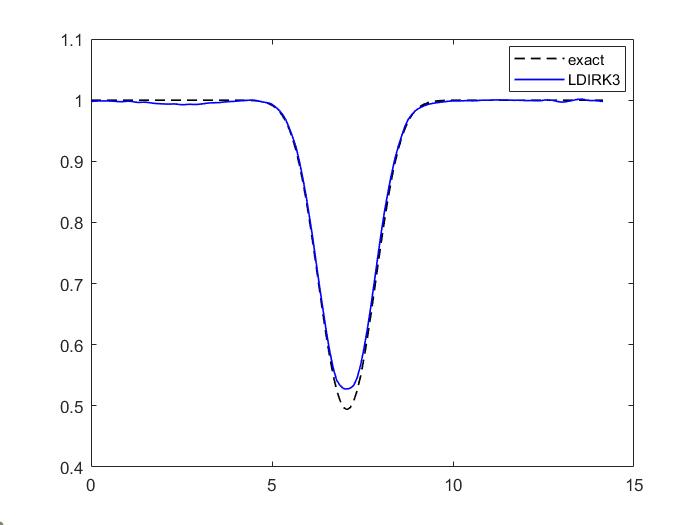}
	\subcaption{L-DIRK3}
	\end{subfigure}
	\caption{Vortex advection problem for 2D Euler equations, density along $x = y$, DIRK3 vs L-DIRK3, $N_x = N_y = 150$, $CFL = 4$, $t = 10$.}
	\label{fig:vortex diagonal}
\end{figure}

 We consider the isentropic vortex advection problem.  The mean flow is $\rho = 1$. $p = 1$, $u = 1$, $v = 1$. The temperature is given by $\displaystyle T = \frac{p}{\rho}$, and the entropy is defined as $S = \frac{p}{\rho^{\gamma}}$. We add the perturbations
\begin{equation*}
	\begin{split}
		&(\delta u, \delta v) = \frac{\epsilon}{2\pi} e^{0.5(1-r^2)}(-\bar{y}, \bar{x}),\\
		& \delta T = -\frac{(\gamma-1)\epsilon^2}{8\gamma\pi^2}e^{1-r^2}, \quad \delta S = 0,
	\end{split}
\end{equation*}
where $(\bar{x}, \bar{y}) = (x-5, y-5)$, $r = \bar{x}^{2}+\bar{y}^{2}$. The vortex strength is taken as $\epsilon = 5$. The exact solution of this problem is simply given by
\[
u(x,y,t) = u(x-t,y-t,0).
\]
The computational domain is set to be $[0, 10]\times[0, 10]$ coupled with periodic boundary conditions in both directions. The vortex is advected along the diagonal direction. Figures \ref{fig:vortex diagonal} compares the solutions of density by DIRK3 and L-DIRK3 along the diagonal $x = y$ after one period of evolution under mesh size $150\times150$ and $CFL = 4$. It is seen that the solution by DIRK3 presents obvious oscillation near the edge of the vortex whereas the result by L-DIRK3 is much less oscillatory. At the center of vortex, the L-DIRK3 is more dissipative since the order of accuracy in time locally drops to first order. Nevertheless, the clipping seems to be acceptable considering the large $CFL$ number applied. Indeed, Table \ref{tab:vortex errors} compares the errors of density of the limited and unlimited schemes. It is seen that applying the time limiter increases the $L^{\infty}$ error due to the stronger dissipation at smooth extrema. However, the $L^1$ and $L^2$ errors are reduced because of the removal of oscillations. 

\begin{table}[h]
	\renewcommand\arraystretch{1.5} 
	\centering
	\begin{tabular}{c|c|c|c}
		\hline
		\hline
		& $L^{\infty}$ err &  $L^{1}$ err & $L^{2}$ err  \\
		\hline
		DIRK3 & $2.39e-2$ & $5.36e-3$ & $6.72e-3$\\
		L-DIRK3& $3.37e-2$ & $2.72e-3$ & $3.95e-3$\\
		\hline
		\hline
	\end{tabular}
	\caption{Vortex advection, density error norms. DIRK3 vs L-DIRK3. $CFL = 4$. 5th order WENO in space, $N_x = N_y = 150$, periodic bc, 1 period of evolution.}
	\label{tab:vortex errors}
\end{table}

\subsection{Summary}

From these numerical results, it is seen that L-DIRK3 is advantageous in the sense that it generates non-oscillatory solutions without damaging the resolution of high order spatial discretization schemes under large CFL numbers. The dimension-by-dimension extension of the scheme to multidimensional problems is validated. Further improvements may be done in several aspects. At first, the solution can still have small amplitude oscillations near discontinuities, which may be removed by improving the design of limiter $\theta^{(k)}_{j}$ and the way of defining $\theta^{(k)}_{\jph}$. For instance, one may construct $\theta^{(k)}_{\jph}$ with different ways of averaging $\theta^{(k)}_{j}$ and $\theta^{(k)}_{j+1}$, such as $\theta^{(k)}_{\jph} = \min(\theta^{(k)}_{j},\theta^{(k)}_{j+1})$. At smooth extrema, the solution has stronger dissipation since the order of accuracy in time is locally dropped to first order. This may be overcome by introducing additional shock detectors so that the time-limiter is switched off in the smooth regions.

\section{Conclusions}

We have reviewed the concepts and motivations of time-limiter schemes. The purpose of the study is to improve the stability of high order implicit time integration methods when the time step is beyond their SSP limits. The
improved stability of time-limiter schemes is realized by taking local convex combination of a first order unconditionally SSP method with a higher order but oscillatory method. Local time-limiters are employed to measure the monotonicity of solution and control the switch between different methods. The idea was applied to construct second order accurate L-Trap and L-DIRK2 schemes \cite{timelimiter}. To better preserve the high resolution of high order space discretizations, one of the main aspects of our work is to apply the idea of limited time integration to improve the stability of the third order accurate, strongly S-stable but non-SSP DIRK3 method. The L-DIRK3 is subsequently constructed. To enhance the flexibility when solving systems of equations, we proposed a new and convenient construction of time-limiters. The new limiter avoids the evaluations of time derivative functions and hence allows arbitrary choice of reference variable with minimal cost. In the previous work by Duraisamy et al \cite{timelimiter} the performance of time-limiter schemes was only tested against one-dimensional conservation laws. To further illustrate the potential applicability in real world applications, another aspect of our study is to extend the discussions of time-limiter schemes to multi-dimensions and convection-diffusion equations, where a large CFL number is desired in order to produce a reasonable time step. The numerical results have shown the advantage of L-DIRK3 in terms of generating high resolution, non-oscillatory solutions under large CFL numbers.

Further investigations are still required at several aspects. At first, the theoretical stability of time-limiter schemes as applied to general nonlinear problems still lacks proofs. The existing framework of SSP methods only applies to time integration methods with constant coefficients. New mathematical tools are needed to enable the nonlinear stability analysis of methods with local, self-adjusting coefficients. Second, with our current constructions some portions of the solution by L-DIRK3 can still exhibit slight oscillations. Indeed, the original DIRK3 scheme is not SSP and it is challenging to completely remove the oscillations near discontinuities. It remains to explore better constructions of time-limiters in order to further remove the oscillations. Moreover, the application of time-limiters may introduce excessive dissipation at extrema and lead to order reduction. Additional shock detecting may be applied to switch off time limiting at smooth extrema.

Overall, our study indicates great potential of L-DIRK3 scheme to be applied in complicated application problems such as DNS. More generally, the framework of time-limiter schemes may be taken as a general and economic approach to improve the stability of an arbitrary DIRK method e.g. TR-BDF2 \cite{TR-BDF2,TR-BDF2analysis,SSPTR-BDF2}.

\bibliographystyle{plain}
\bibliography{references}

\begin{thebibliography}{10}

\bibitem{DIRK}
Roger Alexander.
\newblock Diagonally implicit runge--kutta methods for stiff ode’s.
\newblock {\em SIAM Journal on Numerical Analysis}, 14(6):1006--1021, 1977.

\bibitem{TR-BDF2}
Randolph~E Bank, William~M Coughran, Wolfgang Fichtner, Eric~H Grosse, Donald~J
  Rose, and R~Kent Smith.
\newblock Transient simulation of silicon devices and circuits.
\newblock {\em IEEE Transactions on Computer-Aided Design of Integrated
  Circuits and Systems}, 4(4):436--451, 1985.

\bibitem{SSPTR-BDF2}
Luca Bonaventura and A~Della Rocca.
\newblock Unconditionally strong stability preserving extensions of the tr-bdf2
  method.
\newblock {\em Journal of Scientific Computing}, 70(2):859--895, 2017.

\bibitem{timelimiter}
Karthikeyan Duraisamy, James~D Baeder, and Jian-Guo Liu.
\newblock Concepts and application of time-limiters to high resolution schemes.
\newblock {\em Journal of scientific computing}, 19(1):139--162, 2003.

\bibitem{SSP}
Sigal Gottlieb, Chi-Wang Shu, and Eitan Tadmor.
\newblock Strong stability-preserving high-order time discretization methods.
\newblock {\em SIAM review}, 43(1):89--112, 2001.

\bibitem{TVD}
Ami Harten.
\newblock High resolution schemes for hyperbolic conservation laws.
\newblock {\em Journal of computational physics}, 135(2):260--278, 1997.

\bibitem{UNO}
Ami Harten and Stanley Osher.
\newblock Uniformly high-order accurate nonoscillatory schemes. i.
\newblock In {\em Upwind and High-Resolution Schemes}, pages 187--217.
  Springer, 1997.

\bibitem{monotone}
Amiram Harten, James~M Hyman, Peter~D Lax, and Barbara Keyfitz.
\newblock On finite-difference approximations and entropy conditions for
  shocks.
\newblock {\em Communications on pure and applied mathematics}, 29(3):297--322,
  1976.

\bibitem{TR-BDF2analysis}
ME~Hosea and LF~Shampine.
\newblock Analysis and implementation of tr-bdf2.
\newblock {\em Applied Numerical Mathematics}, 20(1-2):21--37, 1996.

\bibitem{Huynh}
Hung~T Huynh.
\newblock Accurate monotone cubic interpolation.
\newblock {\em SIAM Journal on Numerical Analysis}, 30(1):57--100, 1993.

\bibitem{jia2019}
Feilin Jia, ZJ~Wang, Rathakrishnan Bhaskaran, Umesh Paliath, and Gregory~M
  Laskowski.
\newblock Accuracy, efficiency and scalability of explicit and implicit fr/cpr
  schemes in large eddy simulation.
\newblock {\em Computers \& Fluids}, 195:104316, 2019.

\bibitem{ldirk3scitech}
Jingcheng Lu and James~D Baeder.
\newblock The high resolution l-dirk3 scheme for conservation laws.
\newblock In {\em AIAA SCITECH 2022 Forum}, page 1075, 2022.

\bibitem{martin2006}
M~Pino Mart{\'\i}n and Graham~V Candler.
\newblock A parallel implicit method for the direct numerical simulation of
  wall-bounded compressible turbulence.
\newblock {\em Journal of Computational Physics}, 215(1):153--171, 2006.

\bibitem{S-stable}
A~Prothero and A~Robinson.
\newblock On the stability and accuracy of one-step methods for solving stiff
  systems of ordinary differential equations.
\newblock {\em Mathematics of Computation}, 28(125):145--162, 1974.

\bibitem{Newton}
Thomas Pulliam.
\newblock Time accuracy and the use of implicit methods.
\newblock In {\em 11th Computational Fluid Dynamics Conference}, page 3360,
  1993.

\bibitem{WENO}
Chi-Wang Shu.
\newblock Essentially non-oscillatory and weighted essentially non-oscillatory
  schemes for hyperbolic conservation laws.
\newblock In {\em Advanced numerical approximation of nonlinear hyperbolic
  equations}, pages 325--432. Springer, 1998.

\bibitem{MP}
A~Suresh and HT~Huynh.
\newblock Accurate monotonicity-preserving schemes with runge--kutta time
  stepping.
\newblock {\em Journal of Computational Physics}, 136(1):83--99, 1997.

\bibitem{MUSCL}
Bram Van~Leer.
\newblock On the relation between the upwind-differencing schemes of godunov,
  engquist--osher and roe.
\newblock {\em SIAM Journal on Scientific and statistical Computing},
  5(1):1--20, 1984.

\end{thebibliography}
\end{document}